\documentclass{article}

\usepackage{arxiv}

\usepackage[utf8]{inputenc} 
\usepackage[T1]{fontenc}    
\usepackage{hyperref}       
\usepackage{url}            
\usepackage{booktabs}       
\usepackage{amsfonts}       
\usepackage{nicefrac}       
\usepackage{microtype}      
\usepackage{lipsum}		
\usepackage{graphicx}
\usepackage[numbers]{natbib}
\usepackage{doi}

\usepackage{amssymb}
\usepackage{amsmath}
\usepackage{float}

\usepackage{multirow}
\usepackage{bm}

\newtheorem{proposition}{Proposition}[section]

\newcommand{\diag}{\operatorname{diag}}

\title{Sampling patterns for Zernike-like bases in non-standard geometries}

\author{ \href{https://orcid.org/0009-0004-9623-1766}{\includegraphics[scale=0.06]{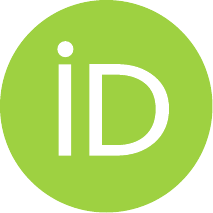}\hspace{1mm}Sergio~Díaz-Elbal} \\
	Department of Mathematics\\
	University and Almería\\
	Almería, Spain \\
	\texttt{sde468@ual.es} \\
	\And
	\href{https://orcid.org/0000-0001-9421-5624}{\includegraphics[scale=0.06]{orcid.pdf}\hspace{1mm}Andrei~Martínez-Finkelshtein} \\
	Department of Mathematics\\
	University and Almería\\
	Almería, Spain \\
	Department of Mathematics\\
	Baylor University\\    
    Waco, Texas, United States \\
	\texttt{andrei@ual.es} \\
	\And
	\href{https://orcid.org/0000-0002-6127-6559}{\includegraphics[scale=0.06]{orcid.pdf}\hspace{1mm}Darío ~Ramos-López} \thanks{Corresponding author.} \\\\
	Department of Mathematics\\
	University and Almería\\
	Almería, Spain \\
	\texttt{dramoslopez@ual.es} \\    
}

\renewcommand{\shorttitle}{Sampling patterns for Zernike-like bases in non-standard geometries}

\hypersetup{
pdftitle={Sampling patterns for Zernike-like bases in non-standard geometries},
pdfsubject={math.NA},
pdfauthor={Sergío~Díaz-Elbal, Andrei~Martínez-Finkelshtein,Darío~Ramos-López},
pdfkeywords={Zernike polynomials, wavefront reconstruction, sampling patterns, interpolation, numerical condition},
}

\begin{document}
\maketitle

\begin{abstract}
	Zernike polynomials are widely used in optics and ophthalmology due to their direct connection to classical optical aberrations. While orthogonal on the unit disk, their application to discrete data or non-circular domains--such as ellipses, annuli, and hexagons--presents challenges in terms of numerical stability and accuracy. In this work, we extend Zernike-like orthogonal functions to these non-standard geometries using diffeomorphic mappings and construct sampling patterns that preserve favorable numerical conditioning. We provide theoretical bounds for the condition numbers of the resulting collocation matrices and validate them through extensive numerical experiments. As a practical application, we demonstrate accurate wavefront interpolation and reconstruction in segmented mirror telescopes composed of hexagonal facets. Our results show that appropriately transferred sampling configurations, especially \textit{Optimal Concentric Sampling} and \textit{Lebesgue points}, allow stable high-order interpolation and effective wavefront modeling in complex optical systems. Moreover, the \textit{Optimal Concentric Samplings} can be computed with an explicit expression, which is a significant advantage in practice.
\end{abstract}


\section{Introduction}
\label{sec:intro}

Zernike polynomials \citep{Zernike1934} are extensively used in both clinical and academic settings in optics and ophthalmology. Their widespread adoption is due to their direct connection to Seidel’s classical optical aberrations \citep{Dai2008}, making them a natural choice for wavefront reconstruction and analysis.

Mathematically, Zernike polynomials form a complete orthogonal basis on the unit disk. However, in practical applications, this orthogonality is often not fully leveraged. One reason for this is that wavefront data is typically sampled at discrete points rather than being continuously defined over the unit disk. Moreover, real-world optical domains are often non-circular, taking the form of polygons, ellipses, annuli, or sectors due to anatomical constraints of the eye or specific optical instrument designs.

In these non-standard domains, standard Zernike polynomials lose their orthogonality and may no longer effectively represent classical optical aberrations. To address this limitation, alternative polynomial bases have been developed \citep{Martinez-Finkelshtein:2009kx, Ares2006, Schneider2009, MR3683670}. When maintaining orthogonality is essential, researchers have proposed various strategies. Some approaches employ the Gram-Schmidt orthogonalization process \citep{Mahajan1994,Mahajan2006,Mahajan2006a}, while others utilize diffeomorphisms to map the unit disk onto non-standard domains \citep{Navarro2014,Ferreira2015}.

Beyond the choice of basis functions, wavefront reconstruction also depends on how the function is fitted to a finite set of measured data points. This process can lead to either an interpolatory function or a discrete best approximation. A crucial factor influencing the accuracy and stability of the reconstruction is the placement of sampling nodes, as their distribution directly affects numerical conditioning and the uniqueness of the solution (unisolvency). Proper selection of these nodes is thus critical for stable and reliable reconstructions. This issue has been extensively studied in the context of the unit disk \citep{Navarro2009,Cuyt2012,RamosLopez2016,Carnicer2014,Meurant2019}. In particular, \cite{RamosLopez2016} introduced a sampling scheme based on concentric circles with optimally determined radii, ensuring an unisolvent set while enhancing numerical stability.

However, determining optimal node configurations for more complex domains remains an open research problem. The primary objective of this work is to discuss node configurations that exhibit favorable numerical stability in various non-standard optical domains, including hexagons, annuli, and ellipses. These configurations are constructed by adapting well-performing sampling patterns from the unit disk to these more complex geometries.

To achieve this, we employ the diffeomorphism technique described in the literature to generate Zernike-like orthogonal bases tailored to each domain. A key focus of our study is the numerical conditioning of these newly constructed function sets when applied to transformed domains. Ensuring well-conditioned basis functions is crucial for achieving accurate and stable wavefront reconstructions in practical applications.

As a specific case study, we analyze hexagonal sampling patterns in a composite domain composed of multiple hexagons. This configuration is particularly relevant in astronomy, where segmented mirror telescopes are constructed using hexagonal tiles \citep{vanKooten2022, Ragland2022}. Our analysis provides insights into the effectiveness of these sampling patterns in such configurations, offering potential improvements for optical system design.

\section{Interpolation with Zernike circle polynomials}
\label{sec:interpolation}

Zernike circle polynomials \cite{Dai2008, Niu2022} are typically expressed in polar coordinates using a double-index notation:
\begin{equation}
\label{eq:zernike_def}
\widetilde{Z}_n^m(\rho,\theta):=\left\{ \begin{array}{lcc} N_n^m R_n^{|m|}(\rho) \cos(m\theta), \quad  & m \geq 0, \\ \\ N_n^m R_n^{|m|}(\rho) \sin(|m|\theta), \quad & m < 0, \end{array} \right.   
\end{equation}
where $n \geq 0$ is the radial order, $|m| \leq n$, and $n-m$ is even. The normalization constant is denoted by $N_n^m$, and the radial component is given by  
\[R_n^{|m|}(\rho):=\sum_{i=0}^{(n-|m|)/2} \frac{(-1)^i(n-i)!}{i! ((n+|m|)/2-i)!((n-|m|)/2-i)!} \rho^{n-2i};\]
it can be expressed in terms of Jacobi orthogonal polynomials or as a hypergeometric function.
 
A single-index notation is also commonly used for Zernike polynomials:
\[
\widetilde{Z}_j(\rho,\theta):=\widetilde{Z}_n^m(\rho,\theta),\] where \[j=\frac{n(n+2)+m}{2}\in \mathbb{N}\cup \{0\},\]
although different conventions exist for indexing.

The functions defined in \eqref{eq:zernike_def} written in Cartesian coordinates $(x,y)$, \footnote{\, Along this text, we, in general, use tilde in the notation of functions in polar coordinates.} 
$$
Z_j(x,y):= \widetilde{Z}_j\left(x(\rho, \theta),y(\rho, \theta)\right)
$$ 
are orthogonal polynomials on $\mathbb{D}$ with respect to the area measure. With the normalization constant $N_n^m=\sqrt{(2(n+1))/(1+\delta_{m,0})}$ in \eqref{eq:zernike_def}, they form an orthonormal set:
\begin{equation}
    \label{orthogonality}
\left\langle Z_j,Z_k\right\rangle = \frac{1}{\pi} \iint_{\mathbb{D}} Z_{j} Z_{k}\, dxdy = \delta_{jk}:=\begin{cases}
    1, & \text{if } j=k, \\
    0, & \text{otherwise.}
\end{cases}
\end{equation}

For a given natural $n\geq0$, the total number of Zernike polynomials with radial order up to $n$ is $N = (n + 1)(n + 2)/2$ \cite{Niu2022}. This corresponds to the dimension of the subspace of bivariate polynomials with total degree less or equal than $n$.  As a consequence of \eqref{orthogonality}, they constitute a basis of that space in $\mathbb{D}$.

In non-trigonometric Fourier analysis, a function $W$ defined on the unit disk $\mathbb{D}$ (such as the wavefront of an optical system with a circular pupil) can be expanded as an infinite series of Zernike polynomials. This representation guarantees convergence in the $L^2$ sense and, under mild regularity conditions, also in a pointwise sense.

In practical applications, however, the series must be truncated to a finite-degree polynomial, yielding the approximation: 
\begin{equation}\label{eq:frontofondatrunc}
      \widetilde{W}(\rho,\theta) \approx \sum_{j=0}^M b_j \widetilde{Z}_{j}(\rho,\theta),
\end{equation}
where the Fourier coefficients $b_j$ have explicit integral formulas over the unit disk $\mathbb{D}$. However, direct computation of these integrals is often infeasible since only a discrete set of sampled values of $ W$ is typically available.

To approximate the coefficients $b_j$, a common approach is least-squares fitting, especially when sampled data is affected by noise. Alternatively, interpolation methods may be preferable for highly accurate data or when modeling specific wavefronts. In this work, we primarily focus on the interpolatory approach.

For a compact domain $\Omega \subset \mathbb{R}^2$ and a fixed degree $n$, an essential problem is selecting optimal sets of interpolation nodes. The case of the unit disk has been widely studied \cite{Navarro2009,Cuyt2012,RamosLopez2016,Carnicer2014,Meurant2019}.

Bivariate polynomial interpolation involves solving a system of linear equations to determine interpolation coefficients. The corresponding collocation matrix, $\bm Z_N:=\left( Z_{i-1}(\bm s_j) \right)_{i,j=1}^N$, where $\bm s_j$ are interpolation nodes, must be well-conditioned for numerical stability. The condition number, $\kappa(\bm Z_N)=\|\bm Z_N\| \|\bm Z_N^{-1}\|$, where $\| \cdot \|$ stands for the $2$-norm, unless stated otherwise, quantifies sensitivity to perturbations in input data. An ill-conditioned matrix amplifies numerical errors, causing significant inaccuracies.

Various node configurations, including spiral, hexagonal, and hexapolar patterns, have been analyzed \cite{Navarro2009}. However, studies such as \cite{RamosLopez2016} indicate that the matrix $ \bm Z_N$ is often ill-conditioned, leading to numerical errors for higher orders.

Fekete points \citep{Bos2010}, which maximize the absolute value of the Vandermonde determinant, are optimal for interpolation but difficult to compute explicitly. \textit{Approximate Fekete points} (AFP), generated using $QR$ factorization-based algorithms \cite{Bos2010}, have been tested in numerical experiments, with available sets limited to order $n=20$ \cite{SommarivaSets}.

\textit{Lebesgue points}, which minimize the Lebesgue constant \cite{Ibrahimoglu2016}, are another popular option, especially for polynomial interpolation, but which also lack an explicit expression. Their approximate positions up to order $n=30$ have been determined numerically in \cite{Meurant2019, SommarivaSets}.

A widely studied family of node configurations based on concentric circles is known as \textit{Bos arrays} \cite{Bos2010}. These configurations are defined as follows: given a maximum polynomial order $n$, radii $r_1 > r_2 > \dots > r_k \geq 0$ are chosen with $k = \lfloor n/2 \rfloor +1$, where $\lfloor \cdot \rfloor$ denotes the integer part (``floor'') function. On each circle of radius $r_j$, $n_j=2n-4j+5$ equally spaced nodes are placed, ensuring a total of $N = (n+1)(n+2)/2$ nodes. The optimal choice of radii remains an open problem.

\textit{Bos arrays} are known to be unisolvent \cite{Bojanov2003}, ensuring the uniqueness of the interpolation solution. Several notable configurations arise from this framework, differing in radius selection. A particular case is described in \cite{Cuyt2012}, the \textit{Cuyt points}, where the radii are expressed in terms of the zeros of Legendre polynomials. 

In \cite{RamosLopez2016}, expressions for radii $r_j$'s were obtained by numerical minimization of the condition number $\kappa(\bm Z_N)$ for different values of $N$. A subsequent least-squares fit of these radii led to the following expression for the numerically optimal values of radii $r_j$:
\begin{equation}
r_j:=r_j(n)=1.1565\xi_{j,n}-0.76535\xi_{j,n}^2+0.60517\xi_{j,n}^3,
\end{equation}
where $\xi_{j,n}$ are zeros of the $(n+1)-$th orthogonal Chebyshev polynomial of the first kind, that is,  
\[\xi_{j,n}=\cos\left(\frac{(2j-1)\pi}{2(n+1)}\right), \quad j=1,\dotsc,k=\left\lfloor\frac{n}{2}\right\rfloor+1.\]

This configuration of nodes, named as \textit{Optimal Concentric Sampling} (OCS), has attracted interest \cite{Meurant2019,Cinzano2020,Dudt2020,Díaz2023} due to its excellent performance in terms of conditioning of the collocation matrix, as well as a growth of the Lebesgue constant, even with very large radial orders (up to $n=30$, which corresponds with 496 polynomials and nodes).

A closely related set of interpolation nodes for Zernike orthogonal polynomials on the disk appears in \cite{Carnicer2014}; we refer to it as the \textit{Carnicer-Godés}, or simply \textit{Carnicer points}. These points, as in the case of OCS sampling and the \textit{Cuyt points}, are part of the family of \textit{Bos arrays}. In this case, the radii are chosen according to the following formula
\begin{equation} \label{eq:Carnicer}
    r_j:=1-\left(\frac{2(j-1)}{n}\right)^a, \quad j=1,\dotsc,k,
\end{equation}
where the exponent $a$ must be determined experimentally and may depend on $n$. A good choice, according to \cite{Carnicer2014}, valid for all values of $n$, is $a=1.46$. It is obtained heuristically from the analysis of the condition number of the collocation matrix $\bm Z_N$ (with the $\infty$-norm) and the Lebesgue constants of the node configuration. 

\begin{table}[tb]
\centering
\scalebox{0.7}{
\begin{tabular}{|c|c|c|c|c|c|c|c|c|c|}
\hline
 & & \multicolumn{8}{c|}{$r_j$} \\
\hline
\multirow{2}{*}{$n=10$} & $OCS$      & 0   & 0.2786 & 0.4972 & 0.6981 & 0.8742 & 0.9818 & - & - \\
\cline{2-10}
& $Carnicer$ & 0   & 0.2780 & 0.5256 & 0.7376 & 0.9046 & 1 & - & - \\
\hline
\multirow{2}{*}{$n=15$} & $OCS$ &  0.1066 & 0.2860 & 0.4385 & 0.5802 & 0.7162 & 0.8398 & 0.9362 & 0.9894\\ 
\cline{2-10}
 & $Carnicer$ & 0.0958 & 0.2780 & 0.4468 & 0.6006 & 0.7376 & 0.8548 & 0.9472 & 1 \\ \hline
\end{tabular}
}
\caption{Radii \(r_j\) of the OCS and \textit{Carnicer points} for orders \(n=10,15\).}
\label{tab:radii}
\end{table}

As shown in Table \ref{tab:radii}, the radii of both node configurations are very similar. However, the OCS radii tend to be more concentrated toward the center of the disk--where the most significant optical information is normally found. On the other hand, \eqref{eq:Carnicer} provides a seemingly simpler formula that, nevertheless, requires determining an optimal constant $a$.

\section{Extension of Zernike polynomials to non-standard geometries}
\label{sec:extensions}

We now examine several extensions of Zernike polynomials--originally defined on the unit disk--to other geometries, such as hexagons, ellipses, and annuli. These domains are particularly relevant due to their frequent appearance in practical applications and the existence of orthogonal bases that can be constructed via diffeomorphic mappings.

There are two main approaches for constructing orthogonal bases on compact planar domains. The first involves generating a basis directly, most commonly using the Gram–Schmidt orthogonalization process \citep{Mahajan1994}. However, this method can lead to significant numerical instability, especially when dealing with high-degree polynomials.

The second approach mitigates these issues by transplanting a known orthogonal basis onto a new domain through a diffeomorphism. This strategy was employed by \cite{Navarro2014}, who used Zernike polynomials on the unit disk as a reference and mapped them onto more general geometries.

Consider a one-to-one mapping 
\begin{equation}
\label{eq:generaldiffeomorphism}
\varphi: \mathbb{D} \rightarrow M \subset \mathbb R^2, \quad (u,v)\rightarrow (x,y)=\varphi(u,v),
\end{equation}
which is a diffeomorphism onto $M=\varphi(\mathbb D)$. By making the corresponding change of variables in \eqref{orthogonality} and considering that $dudv=|J(x,y)|dxdy$, where $J$ is the Jacobian of $\varphi^{-1}$, we obtain that
$$ 
  \frac{1}{\pi} \iint_{M} Z_j(\varphi^{-1}(x,y))Z_k(\varphi^{-1}(x,y))|J_{\varphi^{-1}}(x,y)|\, dxdy=    \delta_{jk}.
$$

In other words, the functions $P_j:=Z_j\circ \varphi^{-1} $ are orthonormal on $M$ with respect to the planar measure $|J(x,y)| dxdy$. 
In general, these functions are no longer polynomials. For that reason, a normalizing factor--a non-vanishing function $q(x,y)$ on $M$--is used, in such a way that the transformed functions
\begin{equation}\label{eq:defzerntransf}
    Q_j(x,y):=q(x,y)P_j( x,y)
\end{equation}
are orthonormal on $M$ with respect to the measure $|J(x,y)|/q^2(x,y)dxdy$.

We now briefly describe the specific transformations for the geometries considered in this study.

\medskip

\noindent \textbf{Hexagonal Zernike basis:} we follow the method of \cite{Ferreira2015} for regular polygons, specializing it for the hexagon inscribed in a unit circle. In the general framework of \eqref{eq:generaldiffeomorphism}, $M$ is a regular $p$-sided polygon. Define $\alpha = \pi/p$ and 

\begin{equation}\label{eq:ralpha}
   R_\alpha(\phi):= \frac{\cos(\alpha)}{\cos (U_\alpha(\theta))}, \quad \text{where} \quad U_\alpha(\theta):=\theta - \left\lfloor \frac{\theta+\alpha}{2\alpha} \right \rfloor 2\alpha, \quad 0\leq \theta \leq 2\pi.
\end{equation}

With the notation \eqref{eq:ralpha}, a diffeomorphism \eqref{eq:generaldiffeomorphism} was constructed in \cite{Ferreira2015}, given in polar coordinates by $\widetilde \varphi(r, \theta) =\left( r R_{\alpha}(\theta), \theta\right)$, with the Jacobian of $\widetilde \varphi^{-1}$  
\[\widetilde{J}(\rho,\theta)=\frac{1}{R_\alpha(\theta)^2}.\]
Correspondingly, two orthonormal families were put forward:
with $q(x,y)=1$ in equation \eqref{eq:defzerntransf}, we obtain the orthonormal family $\{K_j\}_{j\geq0}$, which in polar coordinates is given by
\begin{equation} \label{eq:K}
\widetilde{K}_j(\rho,\theta):=\widetilde{Z}_j\left(\frac{\rho}{R_\alpha(\theta)},\theta\right).
\end{equation}

Alternatively, the choice $q(x,y)=1/R_\alpha(\arctan(y/x))$ yields the family $\{H_j\}_{j\geq0}$, given in polar coordinates by
\begin{equation}\label{eq:H}
    \widetilde{H}_j(\rho,\theta):=\frac{1}{R_\alpha(\theta)}\widetilde{Z}_j\left(\frac{\rho}{R_\alpha(\theta)},\theta\right).
\end{equation}
Clearly, neither of the orthogonal bases \eqref{eq:K}--\eqref{eq:H} 
are polynomial bases.

\medskip

\begin{figure}[H]
   \centering
    \includegraphics[width=8cm]{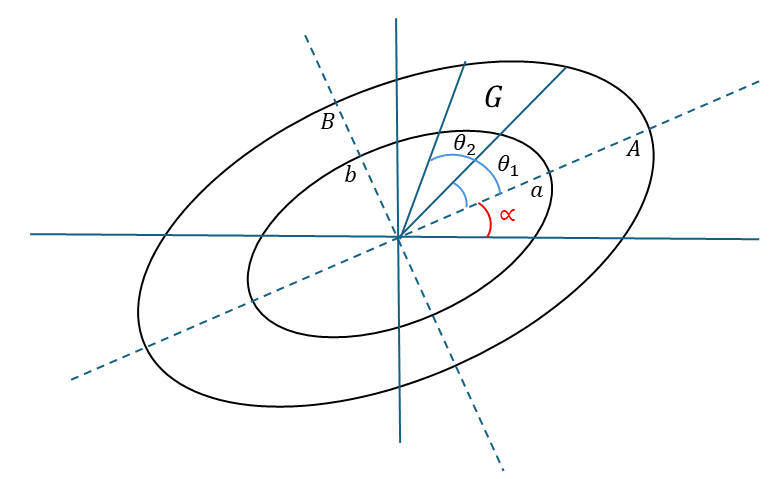}
    \caption{General parametrization of an elliptical aperture $G$.}
    \label{fig:difeomorfismoelipse}
\end{figure}

\noindent \textbf{Elliptical Zernike basis:}
we adopt the general formulation presented in \cite{Navarro2014} for mapping Zernike polynomials onto elliptical sectors, annuli, or a combination of them. An arbitrary elliptical domain $G$ is characterized by multiple parameters (see Figure \ref{fig:difeomorfismoelipse}), providing a general transformation that can be readily adapted to specific geometries. 

For full ellipses aligned with the axes (setting the parameters $a=b=0$, $\theta_1=0,\theta_2=2\pi$, and $h=a/A=0$ in \cite{Navarro2014}), the transformation \eqref{eq:generaldiffeomorphism} boils down to an affine mapping, corresponding to an axes scaling. The Jacobian of this inverse transformation $\varphi^{-1}$ is constant, $J(x,y)=1/AB$, and acts as a normalization factor. Using this technique, we derive the family of polynomials $\{E_j\}_{j\geq0}$, which are orthogonal on the ellipse. Choosing $q(x,y)=1/\sqrt{AB}$ in equation \eqref{eq:defzerntransf}, we obtain, in cartesian coordinates:
\begin{equation} \label{eq:E}
   E_j(x,y):=\frac{1}{\sqrt{AB}} Z_j \left( \frac{x}{A},\frac{y}{B}\right).
\end{equation}

\medskip

\noindent \textbf{Annular Zernike basis:}
using the same approach as in the elliptical case above, the arbitrary elliptical domain $G$ in Figure \ref{fig:difeomorfismoelipse} can be particularized to a circular annulus by setting the parameters $a=b$, $A=B$, $\theta_1=0$, $\theta_2=2\pi$, and $h=a/A$ in \cite{Navarro2014}. Thus, we obtain a transformation as in \eqref{eq:generaldiffeomorphism} with $M$ being the circular annulus. The Jacobian of this inverse transformation $\varphi^{-1}$, expressed in polar coordinates, is given by
\[\widetilde{J}(\rho,\theta)=\frac{\rho-hA}{\rho A^2(1-h)^2}.\]

Following the approach used for the hexagonal domain, we derive two families of orthogonal functions, depending on the choice of $q(x,y)$ in \eqref{eq:defzerntransf}. For $q(x,y) = \sqrt{|\widetilde{J}(\rho,\theta)|}$, one has the orthogonal  basis $\{ O_j \}_{j \geq 0}$, with 
\begin{align} \label{eq:O}
 \widetilde{O}_j(\rho,\theta):=\sqrt{\frac{\rho-hA}{\rho A^2(1-h)^2}} \, {Z}_j\left( \frac{\rho-hA}{A(1-h)}\cos(\theta),\frac{\rho-hA}{A(1-h)} \sin(\theta)\right).
\end{align}
Alternatively, choosing $q(x,y)=1$, one obtains the basis $\{ C_j \}_{j \geq 0}$, with
\begin{equation} \label{eq:C}
    \widetilde{C}_j(\rho,\theta):={Z}_j\left( \frac{\rho-hA}{A(1-h)}\cos(\theta),\frac{\rho-hA}{A(1-h)} \sin(\theta)\right).
\end{equation}
As in the hexagonal case, the orthogonal bases \eqref{eq:O} and \eqref{eq:C} are not polynomial.

\section{Good sampling patterns in non-standard geometries}
\label{sec:samplings}
 
The main idea is to transplant the sampling points from the disk to the target domain $M$ in \eqref{eq:generaldiffeomorphism} employing the same diffeomorphism $\varphi$ used to construct the orthogonal functions. More precisely, if 
$$
\mathcal S=\{\bm s_1, \bm s_2,\dotsc,\bm s_N\} \subset \mathbb{D} 
$$ 
is a node configuration on the unit disk, where $N$ is the number of Zernike polynomials of order $\le n$, then 
\begin{equation}\label{eq:transformhex}
    \varphi(\mathcal S)=\{\varphi(\bm s_1),\varphi(\bm s_2),\dotsc,\varphi(\bm s_N)\}
\end{equation}
is a node configuration on $M$. 

With this choice of nodes, the collocation matrices of functions $P_j=Z_j\circ \varphi^{-1} $ preserve the condition number of the original nodes in the disk for Zernike polynomials. This result, whose proof is immediate, can be expressed as:

\begin{proposition}
    \label{prop:generalconditioning}
    Consider a diffeomorphism $\varphi: \mathbb{D} \rightarrow M \subset \mathbb{R}^2$, and the functions $P_j=Z_j\circ \varphi^{-1}$ defined on $M$, where $Z_j$ are the regular Zernike polynomials on the disk. Then, the collocation matrix of functions $P_j$ at the transferred nodes $\varphi(\mathcal S)$ verifies
    \[\left( P_{i-1}(\varphi(\bm s_j) )\right)_{i,j=1}^N= \left( Z_{i-1}(\bm s_j) \right)_{i,j=1}^N = \bm Z_N. \] 
    In particular, the condition number $\kappa(\cdot)$ remains invariant for any underlying matrix $p$-norm. The same is true for functions $Q_j$ from \eqref{eq:defzerntransf} if the factor $q$ is constant. 
\end{proposition}

For functions $Q_j$ from \eqref{eq:defzerntransf} with a non-constant factor $q$,
the collocation matrix $\left(Q_{i-1}(\varphi(\bm s_j) )\right)_{i,j=1}^N $ is equal to $\bm Z_N$, multiplied on the right by a diagonal matrix made of evaluations of $q$ at $\varphi(\mathcal S) $, operation that is known as diagonal preconditioning. 
In this case (which is the case of functions $ H_j $ and $\widetilde O_j $, defined in \eqref{eq:H} and \eqref{eq:O}, respectively), the condition number of their collocation matrix may vary. This situation will be addressed below in more detail. 

\medskip

\begin{figure}[H]
    \centering
    \includegraphics[width=0.47 \linewidth]{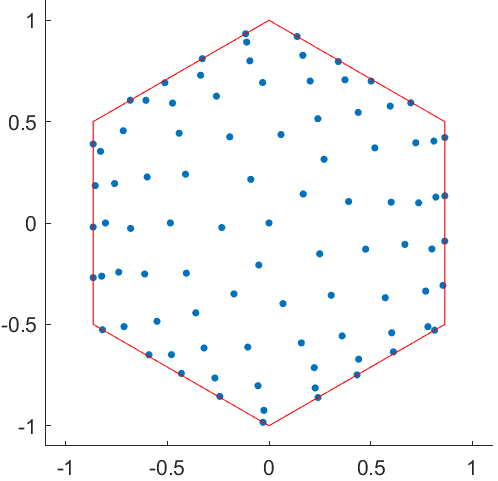}
    \includegraphics[width=0.47 \linewidth]{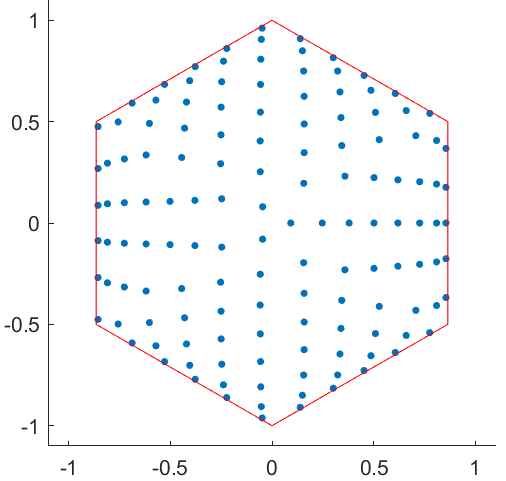}
    \caption{\textit{Lebesgue points} of degree $n=12$ (left) and OCS of degree $n=15$ (right) transferred to the hexagon.}
    \label{fig:hexlebesgue}
\end{figure}

\noindent \textbf{Configurations on the hexagon:} the diffeomorphism of the unit disk $\mathbb{D}$ onto its inscribed hexagon of side $1$ is given by \eqref{eq:ralpha} with $\alpha=\pi/6$. Figure \ref{fig:hexlebesgue} depicts the transformed nodes \eqref{eq:transformhex} corresponding to the original configuration $\mathcal S$ in the disk for the OCS (left) and the \textit{Lebesgue points} (right). 

For functions $K_j$ in \eqref{eq:K}, the collocation matrix over $\varphi(\mathcal S)$ satisfies 
$$
\bm K_N:=\left( K_{i-1}(\varphi(\bm s_j) )\right)_{i,j=1}^N=\bm Z_N, 
$$ 
as in Proposition \ref{prop:generalconditioning}, so that its numerical conditioning is preserved. 

The case of functions $H_j$, where Proposition \ref{prop:generalconditioning} does not apply, must be analyzed separately.

\begin{proposition}   \label{prop:condhex} 
    Consider the collocation matrix \[ \bm H_N:=\left(H_{i-1}(\varphi(\bm s_j)\right)_{i,j=1}^N,\] at the transferred nodes $\varphi(\mathcal S)$. Let $\theta_j:=\arctan(y_j/x_j)$, where $\bm s_j = (x_j,y_j)$, and consider the diagonal matrix \[\bm D_N:=\diag\left(\frac{1}{R_\alpha(\theta_j)}\right), \quad j=1,\dotsc,N.\] Then, 
    \[
    \kappa_p(\bm H_N)\leq 
    \frac{2\sqrt{3}} {3} \, \kappa_p(\bm Z_N) \approx 1.15 \, \kappa_p(\bm Z_N), 
    \] 
    where $\kappa_p$ is any matrix condition number induced by a p-norm. 
\end{proposition}
\textbf{Proof:}
as it was noticed above, we can write $\bm H_N$ as 
\begin{align*}
    \bm H_N &= \underbrace{\begin{pmatrix}
        {Z_0(\bm s_1)}&  {Z_0(\bm s_2)} & \cdots &  {Z_0(\bm s_N)} \\ 
        {Z_1(\bm s_1)}& {Z_1(\bm s_2)}& \cdots &  {Z_1(\bm s_N)}\\
        \vdots & \vdots & \ddots & \vdots \\  {Z_N(\bm s_1)}& {Z_N(\bm s_2)} & \cdots &  {Z_N(\bm s_N)}
    \end{pmatrix}}_{\bm Z_N} \underbrace{\begin{pmatrix}
       \frac{1}{R_\alpha(\theta_1))} &   &  &   \\ 
         &  \frac{1}{R_\alpha(\theta_2))}  &  &   \\ 
        & & \ddots &  \\  &  & & &  \frac{1}{R_\alpha(\theta_N))}
    \end{pmatrix}}_{\bm D_N}.
\end{align*}
Function $R_\alpha(\cdot)$ in \eqref{eq:ralpha} verifies that $\cos(\alpha) \leq R_\alpha(\cdot)\leq 1$, and then, with $\alpha=\pi/6$,
\[
1 \leq \frac{1}{R_{\pi/6}(\cdot)}\leq \frac{1}{ \cos(\pi/6)}=\frac{2\sqrt{3}}{3}.
\] For a diagonal matrix, any matrix $p$-norm $\|\cdot\|_p$ is equal to the largest absolute value of elements in the diagonal. Using this fact we can conclude that
\[
\|\bm D_N\| \cdot \|\bm D_N^{-1}\|=\frac{\max_{j}(1/R_{\pi/6}(\theta_j))}{\min_{k}(1/R_{\pi/6}(\theta_k))}\leq \frac{2\sqrt{3}}{3}.\]
In particular, for the condition number $\kappa_p(\cdot)$, induced by the $p$-norm $\|\cdot\|_p$, we have
\begin{align*}
    \kappa_p(\bm H_N)&=\|\bm H_N\|_p \cdot\|\bm H_N^{-1}\|_p = \|\bm Z_N \bm D_N \|_p \cdot \|(\bm Z_N \bm D_N)^{-1}\|_p \\ &\leq \|\bm Z_N\|_p \cdot\|\bm Z_N^{-1}\|_p \cdot \|\bm D_N\|_p \cdot \|\bm D_N^{-1}\|_p \\ &\leq 
    \frac{2\sqrt{3}}{3}\, \kappa_p(\bm Z_N).
\end{align*}
$\hfill\blacksquare$

Proposition \ref{prop:condhex} shows that the condition number of $\bm H_N$ is in the same order of magnitude as $\bm Z_N$, which guarantees the numerical stability of the translated samplings for both orthogonal bases, $\{ K_j\}_{j\ge 0}$ and $\{ H_j\}_{j\ge 0}$, in the regular hexagon.

\medskip

\noindent \textbf{Configurations on the ellipse:} the transformation that sends the unit disk $\mathbb{D}$ to the ellipse of major axis $A$ and minor axis $B$ is an affine transformation that can be expressed in Cartesian coordinates as 
\begin{equation}
  (x,y)=\varphi(u,v)= (Au,Bv), \quad u^2+v^2\leq 1.
\end{equation}

\begin{figure}[H]
    \centering
    \includegraphics[width=0.47 \linewidth]{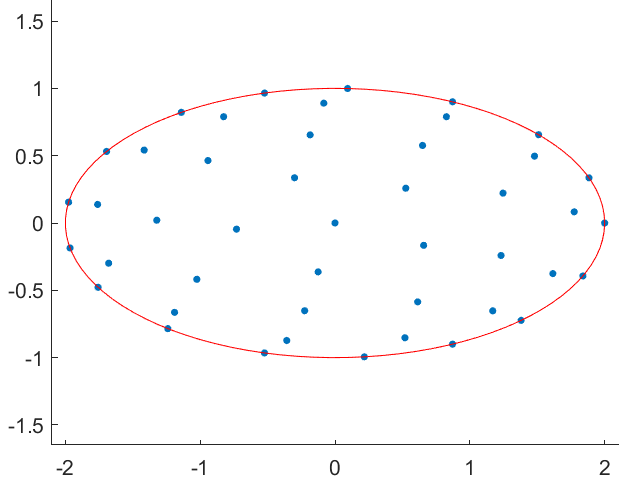}
    \includegraphics[width=0.47 \linewidth]{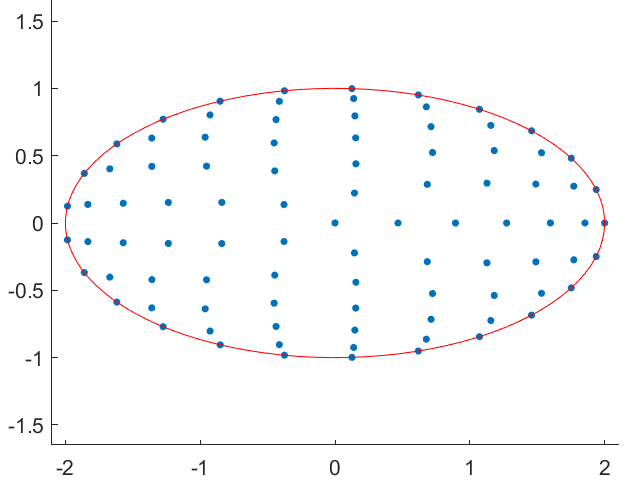}
    \caption{\textit{Fekete points} of degree $n=8$ and  \textit{Carnicer points} of degree $n=12$ transferred to the ellipse.}
    \label{fig:ellfekete}
\end{figure}

As an illustration of this transformation, in Figure \ref{fig:ellfekete} we depict the transferred \textit{Fekete} and \textit{Carnicer points}.

For the orthogonal basis $\{ E_j \}_{j\geq 0} $ on the ellipse, as a consequence of Proposition \ref{prop:generalconditioning}, we have 
\[\ \bm E_N:= \left(E_{i-1}(\varphi(\bm s_j) )\right)_{i,j=1}^N=\bm Z_N, \] and thus the numerical conditioning of the nodes in the disk is preserved.

\medskip

\noindent \textbf{Configurations on the annulus:}
the transformation between the unit disk $\mathbb{D}$ and the annulus of inner radius $a$ and outer radius $A$, has the expression
\begin{equation}
  (x,y)=\varphi(u,v)= (A(h+(1-h)\rho)\cos(\theta),A(h+(1-h)\rho)\sin(\theta)), \quad u^2+v^2\leq 1,
\end{equation}
where $h=a/A$ and $(\rho, \theta)$ are the polar coordinates of $(x, y)$.

\begin{figure}[H]
    \centering
    \includegraphics[width=0.47 \linewidth]{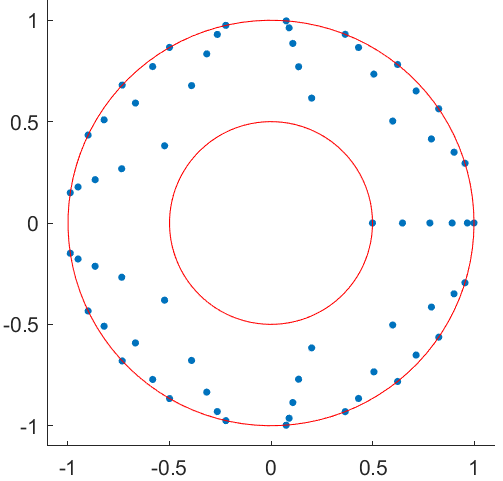}
     \includegraphics[width=0.47 \linewidth]{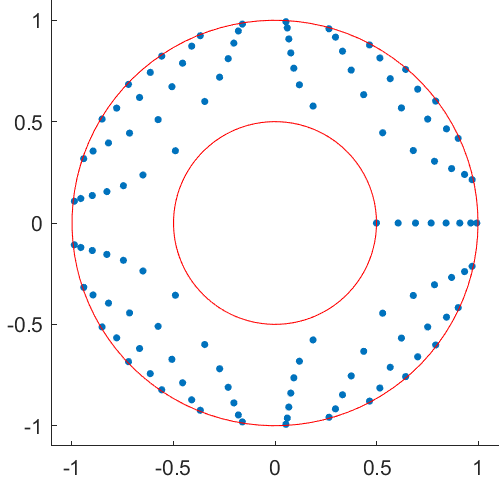}
    \caption{\textit{Cuyt points} of degree $n=10$ and OCS of degree $n=14$, transferred to the annulus.}
    \label{fig:corcuyt}
\end{figure}

As a sample, Figure \ref{fig:corcuyt} shows the \textit{Cuyt points} and the OCS configuration transferred to the annulus.

Similarly to other cases discussed above, for the orthogonal basis $\{C_j\}_{j\geq0}$ given in \eqref{eq:C}, Proposition \ref{prop:generalconditioning} states that
$$
\bm C_N := \left(C_{j-1}(\varphi(\bm s_k)\right)_{j,k=1}^N=\bm Z_N,
$$
and thus its condition number over the nodes $\varphi(\mathcal{S})$ is also invariant. 

In contrast, for the other orthogonal family in the annulus, $\{O_j\}_{j\geq0}$ given in \eqref{eq:O},  
$$
\bm O_N := \left(O_{j-1}(\varphi(\bm s_k)\right)_{j,k=1}^N=\left(\sqrt{\frac{\rho-hA}{\rho A^2(1-h)^2}}Z_{j-1}(\bm s_k)\right)_{j,k=1}^N,
$$
the analysis of $\kappa(\bm O_N)$ is not straightforward. On the one hand, the function  
\begin{equation}
    \label{eq:annulus_q}
    q(x,y) = \sqrt{|\widetilde{J}(\rho,\theta)|} = \sqrt{\frac{\rho-hA}{\rho A^2(1-h)^2}}
\end{equation} 
vanishes when $\rho=a$ (the inner circle of the annulus), which is mapped to the origin in the disk. In this case, the collocation matrix becomes singular and the interpolation is ill-defined. For $\rho \rightarrow a$, the condition number $\kappa(\bm O_N) \rightarrow +\infty$. In order to overcome this problem, we have experimentally tested the effect of modifying node configurations for the nodes corresponding to $\rho_k=a$. This is analyzed in Section \ref{sec:experiments}.

On the other hand, when the inner radius $a$ tends to $A$, the function \eqref{eq:annulus_q} also tends to infinity. Scenarios in which $a$ is very close to $A$ could cause numerical instability. Fortunately, this situation rarely occurs in real applications.

\section{Numerical experiments}
\label{sec:experiments}

We conducted a series of experiments to evaluate the numerical conditioning of the interpolation problem using our proposed sampling schemes for each type of domain, comparing them with other relevant sets of nodes. For each interpolation of order $n$, we selected $N =  (n+1)(n+2)/2$ nodes, corresponding to the transformed configuration $\varphi(\mathcal{S})$ within the target domain $M$. At each node, we evaluated the associated orthogonal basis functions--constructed in Section \ref{sec:extensions}--up to degree $n$, and assembled the resulting $N \times N$ collocation matrix. We then computed its condition number, $\kappa_2(\cdot)$ induced by the matrix $2$-norm.

This process was repeated for interpolation orders up to $n = 30$ (yielding 496 nodes and basis functions), except in the case of \textit{Fekete points}, for which our dataset is limited to orders up to $n = 20$.

\medskip

\noindent \textbf{Configurations on the unit disk and the ellipse:}
as an introduction--and to provide a baseline for comparison with subsequent results--we first analyzed the condition number of the collocation matrix for various sampling patterns on the unit disk, using the Zernike polynomials $\{Z_j\}_{j \geq 0}$. The results are summarized in Table \ref{tab:discocond}, which generally indicate favorable conditioning of the collocation matrix $\bm{Z}_N$ across all sampling schemes, as reflected by the condition number $\kappa_2$.

\begin{table}[tb]
\begin{center}
\scalebox{0.7}{
   $\begin{array}{|c|c|c|c|c|c|}\hline  n & Fekete & Lebesgue & Cuyt & Carnicer & OCS \\ \hline 1 & 1.5399 & 1.4142 & 1.4142 & 1.4142 & 1.0894\\ 2 & 2.7439 & 2.7324 & 2.7324 & 2.7324 & 1.3050\\ 3 & 3.2536 & 3.1225 & 3.1566 & 3.1581 & 1.7631\\ 4 & 3.4446 & 3.3015 & 3.3225 & 3.2647 & 2.0453\\ 5 & 4.6120 & 3.9513 & 4.3097 & 4.2150 & 2.4867\\ 6 & 4.9542 & 4.0784 & 4.6605 & 4.4578 & 2.7353\\ 7 & 5.6333 & 4.9161 & 5.3029 & 5.1572 & 3.2308\\ 8 & 6.3214 & 5.1008 & 6.1798 & 5.5954 & 3.4889\\ 9 & 8.3085 & 5.3493 & 7.3072 & 6.2337 & 4.0410\\ 10 & 8.2548 & 6.3572 & 8.8531 & 6.9373 & 4.3396\\ 11 & 11.9448 & 6.7369 & 10.8302 & 7.7707 & 4.9642\\ 12 & 17.8157 & 7.0520 & 13.5731 & 8.7855 & 5.3384\\ 13 & 18.5226 & 8.0623 & 17.1568 & 9.9757 & 6.0638\\ 14 & 12.6567 & 8.4132 & 22.0487 & 11.4031 & 6.5551\\ 15 & 14.0587 & 8.4660 & 28.4848 & 13.0202 & 7.4148\\ 16 & 17.1278 & 9.4538 & 37.2457 & 14.9912 & 8.0713\\ 17 & 31.5826 & 9.8407 & 48.8657 & 17.2379 & 9.1269\\ 18 & 33.1363 & 10.0008 & 64.6647 & 19.9636 & 10.0257\\ 19 & 29.6446 & 10.9025 & 85.7326 & 23.1097 & 11.3638\\ 20 & 46.7642 & 10.4536 & 114.3812 & 26.9060 & 12.6065\\ 21 & - & 11.2767 & 152.7292 & 31.3207 & 14.3577\\ 22 & - & 12.0337 & 204.8834 & 36.7983 & 16.1049\\ 23 & - & 12.9792 & 274.9192 & 44.4349 & 18.4636\\ 24 & - & 13.0622 & 370.2150 & 54.7430 & 20.9946\\ 25 & - & 13.8224 & 498.5231 & 67.0134 & 24.2573\\ 26 & - & 14.4974 & 673.2165 & 83.3803 & 28.7151\\ 27 & - & 16.2366 & 908.9529 & 102.8985 & 34.0948\\ 28 & - & 15.2266 & 1.2301e+03 & 128.9610 & 40.7343\\ 29 & - & 16.4605 & 1.6643e+03 & 160.1399 & 48.8196\\ 30 & - & 16.3907 & 2.2562e+03 & 201.7801 & 58.7650 \\ \hline   \end{array}$
   }
   \caption{Condition numbers $\kappa$ of the collocation matrix $\bm Z_N=\left( Z_{i-1}(\bm s_j) \right)_{i,j=1}^N$ in each of the samplings for the unit disk.}
   \label{tab:discocond}
   \end{center}
\end{table}

The results in Table \ref{tab:discocond} also provide insight into the behavior on the ellipse and the orthogonal polynomials $\{E_j\}_{j \geq 0}$ constructed in \eqref{eq:E}. Proposition \ref{prop:generalconditioning} ensures that the numerical conditioning in this case is identical to that observed on the unit disk.

\medskip

\noindent \textbf{Configurations on the hexagon:}
Proposition \ref{prop:generalconditioning} guarantees that the conditioning of the samples employing the functions $K_j$, defined in \eqref{eq:K}, is identical to the corresponding one in the unit disk (see Table \ref{tab:discocond}). 
For the case of $\{ H_j\}_{j\geq0}$ given in \eqref{eq:H} (which differs from $K_j$ in the factor $1/R_{\alpha}(\cdot)$), the numerical conditioning of the collocation matrix is virtually unaffected, as predicted by Proposition \ref{prop:condhex}. The results of the experiments are shown in Table \ref{tab:discohex}. 

\begin{table}[tb] \begin{center} \scalebox{0.7}{ $\begin{array}{|c|c|c|c|c|c|} \hline n & Fekete & Lebesgue & Cuyt & Carnicer & OCS \\ \hline  1 & 1.5424 & 1.4142 & 1.4142 & 1.4142 & 1.0894\\ 2 & 2.6644 & 2.6381 & 2.6381 & 2.6381 & 1.3022\\ 3 & 3.4047 & 3.0600 & 3.0648 & 3.0662 & 1.7487\\ 4 & 3.4646 & 3.3446 & 3.3502 & 3.3665 & 2.0837\\ 5 & 4.7611 & 3.8632 & 4.2325 & 4.1088 & 2.5098\\ 6 & 5.1265 & 4.1627 & 4.7841 & 4.5872 & 2.7742\\ 7 & 5.6318 & 4.9381 & 5.3100 & 5.1008 & 3.2715\\ 8 & 6.4440 & 5.3006 & 6.3034 & 5.7390 & 3.5458\\ 9 & 8.4318 & 5.5113 & 7.4074 & 6.2807 & 4.1172\\ 10 & 8.2043 & 6.6362 & 8.9542 & 7.0857 & 4.4395\\ 11 & 12.2404 & 6.7054 & 11.0277 & 7.8981 & 5.0903\\ 12 & 18.6916 & 7.1561 & 13.8663 & 8.9648 & 5.4859\\ 13 & 18.9939 & 8.2517 & 17.6401 & 10.2007 & 6.2447\\ 14 & 12.3993 & 8.6073 & 22.7265 & 11.6907 & 6.7554\\ 15 & 14.5990 & 8.9438 & 29.5367 & 13.4182 & 7.6654\\ 16 & 17.8455 & 9.2566 & 38.7135 & 15.4657 & 8.3408\\ 17 & 30.3842 & 10.1587 & 51.0616 & 17.8722 & 9.4625\\ 18 & 34.4056 & 10.4381 & 67.7225 & 20.7151 & 10.3709\\ 19 & 30.0170 & 11.5052 & 90.2043 & 24.0800 & 11.7973\\ 20 & 48.9868 & 10.9410 & 120.5793 & 28.0594 & 13.0385\\ 21 & - & 11.8041 & 161.6430 & 32.7823 & 14.9107\\ 22 & - & 12.1411 & 217.2085 & 38.3790 & 16.6391\\ 23 & - & 13.1619 & 292.4459 & 45.1537 & 19.1885\\ 24 & - & 13.0789 & 394.3992 & 55.5390 & 22.0302\\ 25 & - & 14.3365 & 532.6510 & 68.0373 & 25.8231\\ 26 & - & 15.0311 & 720.2371 & 84.7075 & 30.5322\\ 27 & - & 16.7740 & 974.9349 & 104.5959 & 36.3070\\ 28 & - & 15.3181 & 1.3209e\!+\!03 & 131.1732 & 43.3877\\ 29 & - & 16.1864 & 1.7912e\!+\!03 & 162.9679 & 52.0749\\ 30 & - & 16.4365 & 2.4308e\!+\!03 & 205.4732 & 62.7200 \\ \hline \end{array}$ } \caption{Condition numbers $\kappa_2$ of the collocation matrices $\bm H_N=\left(H_{i-1}(\bm s_j)\right)_{i,j=1}^N$ for each node configuration mapped to the hexagon.} \label{tab:discohex} \end{center} \end{table}

\medskip

\noindent \textbf{Configurations on the annulus:}
again, Proposition \ref{prop:generalconditioning} guarantees that the conditioning of the samplings with the functions $C_j$, given in \eqref{eq:C}, is identical to the unit disk, see Table \ref{tab:discocond}. 

\begin{table}[tb] \begin{center} \scalebox{0.7}{ $\begin{array}{|c|c|c|c|c|c|} \hline n & Fekete & Lebesgue & Cuyt & Carnicer & OCS \\ \hline 1 & 1.5399 & 1.4142 & 1.4142 & 1.4142 & 1.0894\\ 2 & 5.6530 & 13.0693 & 13.0693 & 13.0693 & 6.2580\\ 3 & 3.8305 & 3.8305 & 3.9425 & 3.9457 & 2.3146\\ 4 & 8.5933 & 9.2258 & 16.3100 & 16.2516 & 8.9572\\ 5 & 6.0460 & 5.3222 & 6.1427 & 6.1439 & 3.6421\\ 6 & 11.5300 & 17.8841 & 19.2637 & 19.3219 & 11.2881\\ 7 & 7.9657 & 6.6345 & 8.2301 & 8.3816 & 5.0701\\ 8 & 14.7927 & 20.3604 & 21.9942 & 22.5124 & 13.4426\\ 9 & 10.0415 & 7.6209 & 10.3285 & 10.8781 & 6.6287\\ 10 & 19.0184 & 23.1542 & 24.5657 & 26.0831 & 15.5836\\ 11 & 13.0292 & 11.2410 & 12.5101 & 13.8407 & 8.3854\\ 12 & 23.2890 & 26.2121 & 27.0232 & 30.4782 & 17.8852\\ 13 & 21.0524 & 13.3994 & 17.7900 & 17.6501 & 10.4483\\ 14 & 26.4920 & 18.5879 & 29.3992 & 36.4852 & 20.5647\\ 15 & 19.8817 & 15.6447 & 29.2002 & 22.9999 & 12.9595\\ 16 & 28.2465 & 29.0063 & 38.0386 & 45.4592 & 23.8555\\ 17 & 35.4970 & 14.8135 & 49.7670 & 31.0772 & 16.1227\\ 18 & 43.0206 & 33.1381 & 65.7066 & 59.6068 & 28.1236\\ 19 & 44.5684 & 18.4383 & 86.9626 & 43.8193 & 20.2433\\ 20 & 51.9107 & 34.7701 & 115.8508 & 82.3726 & 33.8176\\ 21 & - & 17.7106 & 154.5142 & 64.3283 & 25.7458\\ 22 & - & 36.8589 & 207.0704 & 119.1053 & 41.5837\\ 23 & - & 24.2163 & 277.6323 & 97.5817 & 33.2713\\ 24 & - & 40.5045 & 373.6002 & 178.2601 & 52.3668\\ 25 & - & 24.9974 & 502.7867 & 151.6353 & 43.7868\\ 26 & - & 41.1193 & 678.6072 & 273.4479 & 67.5685\\ 27 & - & 30.6723 & 915.8179 & 239.6360 & 58.7493\\ 28 & - & 45.8095 & 1.2389e\!+\!03 & 426.8275 & 89.2716\\ 29 & - & 32.4057 & 1.6756e\!+\!03 & 383.1633 & 80.3447\\ 30 & - & 48.6714 & 2.2707e\!+\!03 & 674.6989 & 120.5633 \\ \hline \end{array}$ } \caption{Condition numbers $\kappa_2$ of the collocation matrices $\bm O_N=\left(O_{i-1}(\bm s_j)\right)_{i,j=1}^N$ for each node configuration mapped to the annulus of inner radius $a=0.5$ and outer radius $A=1$.}\label{tab:tablacorona} \end{center} \end{table}

For the family $\{O_j\}_{j \geq 0}$ introduced in equation~\eqref{eq:O}, the normalization prefactor diverges as $\rho \to 0$, which negatively impacts the numerical conditioning of the resulting collocation matrix--particularly when the original sampling includes a node at the origin. In the case of the annulus, the origin $(0,0)$ is mapped to the point $(a,0)$, which lies on the inner boundary of the domain. A common remedy is to perturb the problematic node slightly, replacing it with $(a + \varepsilon, 0)$ for a small $\varepsilon > 0$. The numerical results, obtained with parameters $a = 0.5$, $A = 1$, and $\varepsilon = 0.01$, are presented in Table~\ref{tab:tablacorona}.
 
As observed, the condition number reaches $10^3$ for \textit{Cuyt points} and approximately an order of magnitude less for OCS, with very similar behavior between the two bases (see also Table~\ref{tab:discocond}). These results are generally favorable, indicating that the sampling configurations--including the small modification at the inner boundary--are practically sound. Moreover, we tested various combinations of inner and outer radii $a$ and $A$, and encountered significant numerical issues only when the two radii are very close, as discussed in Section~\ref{sec:samplings}.

\section{Wavefront fitting of atmospheric turbulence}
\label{sec:wavefrontfit}

\begin{figure}[tb]
    \centering
    \includegraphics[width=0.47 \linewidth]{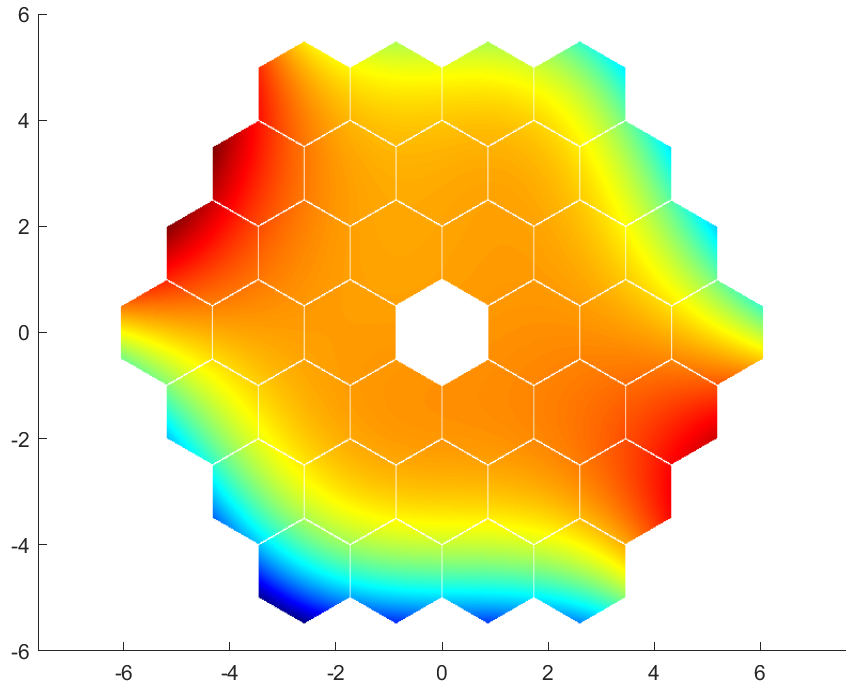}
    \includegraphics[width=0.47 \linewidth]{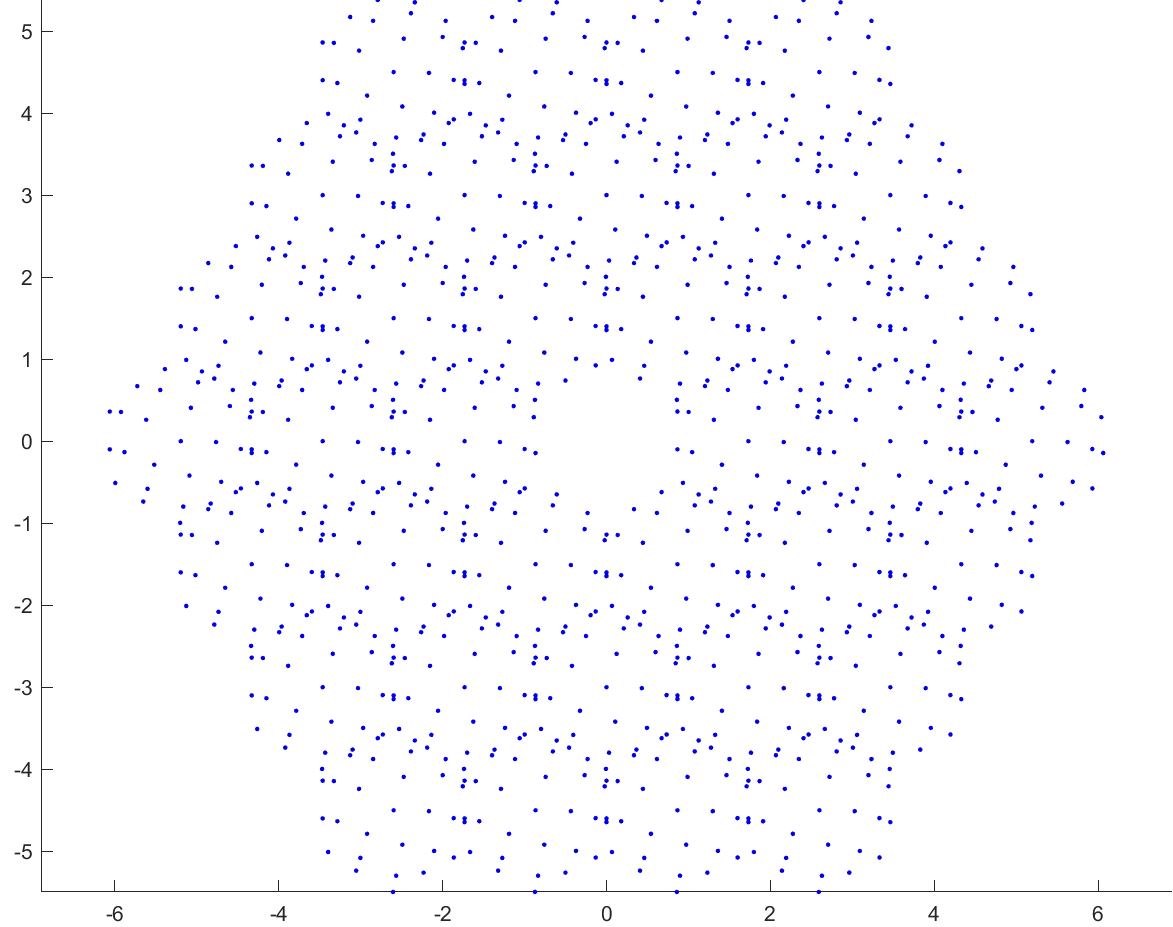}
    \caption{Random atmospheric turbulence following Kolmogorov model (left) and configuration of nodes from \textit{Lebesgue points} of order $n=6$ in each of the 36 hexagons. }
    \label{fig:smtrandom}
\end{figure}

One of the goals of this work is to apply the proposed sampling patterns and the adapted Zernike orthogonal bases to the problems of interpolation and wavefront reconstruction in optical systems, modeled as surfaces in $\mathbb{R}^3$. In \cite{Navarro2014}, wavefront reconstruction was performed over an aperture $\mathcal{H}$ composed of 36 hexagonal sub-apertures, arranged as shown in Figure~\ref{fig:smtrandom}. In our study, we reconstruct a wavefront 
\begin{equation}
\label{eq:kolmo}
     \widetilde{f}(\rho,\theta):=\sum_{j=1}^{14} a_j \widetilde{Z}_{j-1}(\rho,\theta), \quad \rho \leq 6,
\end{equation} 
that is randomly generated from a set of Zernike coefficients $a_j$. Function $\widetilde{f}$ in \eqref{eq:kolmo} can be seen directly as a function defined in the composed domain $\mathcal{H}$. Thus, our problem is to interpolate $\widetilde{f}$ on the domain $\mathcal{H}$.

The coefficients $a_j$ in \eqref{eq:kolmo} are specifically selected according to the Kolmogorov model, in order to generate realistic and random atmospheric turbulence effects \cite{Roddier1990}. The reconstruction is carried out by interpolating the wavefront at a set of nodes within each hexagonal sub-aperture, where the function values are assumed to be known.

In \cite{Navarro2014}, the reconstruction is performed using discrete least squares and an integral approach. In contrast, our method relies on interpolation via \textit{critical sampling} for each order \( n = 1, \dotsc, 30 \). This means that the number of nodes matches the dimension of the polynomial subspace used for approximation.

Therefore, for each order $n$ we select $N=(n+1)(n+2)/2$ nodes arranged according to a chosen sampling configuration within each of the 36 hexagons. The final reconstructed surface is obtained by independently interpolating over each hexagon, resulting in a local--or zonal--interpolation across the entire domain.

To adapt the functions $K_j$ and $H_j$ from Section~\ref{sec:extensions} to the position of each hexagon, we apply a translation to the basis functions. Specifically, for each hexagon ($ k = 1, \dotsc, 36$), with center at $(x_k, y_k)$, we define the shifted functions as
$$
Q_j^{(k)}(x, y) := Q_j(x - x_k, y - y_k) \, \chi_k(x, y),
$$
where $Q_j$ denotes either $K_j$ or $H_j$, and $\chi_k$ is the characteristic function of the $k$-th hexagon.
This transformation effectively translates each orthogonal basis function to the corresponding hexagon in the domain, ensuring that the interpolation is localized. The layout and construction of the system of 36 hexagons are described in detail in \cite{Navarro2009}. In this way,  the interpolating function is 
\[\widetilde{f}^{approx}(\rho,\theta):=\sum_{k=1}^{36} \sum_{j=1}^N c_j^{k} \widetilde{Q}_{j-1}^{(k)}(\rho_k,\theta_k). \]
Furthermore, we sample the wavefront on a non-uniform grid. To construct this grid, we transfer the \textit{Lebesgue points} from \cite{Meurant2019} and the OCS points from \cite{RamosLopez2016} onto the central hexagon, and then replicate this configuration across each of the 36 hexagons that form the aperture. The resulting reconstruction is shown in Figure~\ref{fig:smtrandom}.

For comparison, we also perform the experiment using a random sampling strategy within the hexagon. We initially generate 1000 random points and apply a thinning algorithm to reduce the set to the number of nodes required for critical sampling at each interpolation order. Additionally, we include a spiral sampling scheme from \cite{Navarro2009}--originally defined on the disk--which we adapt to the hexagonal geometry through a suitable transformation. 

\begin{figure}[tb]
    \centering
    \includegraphics[width=6.7cm]{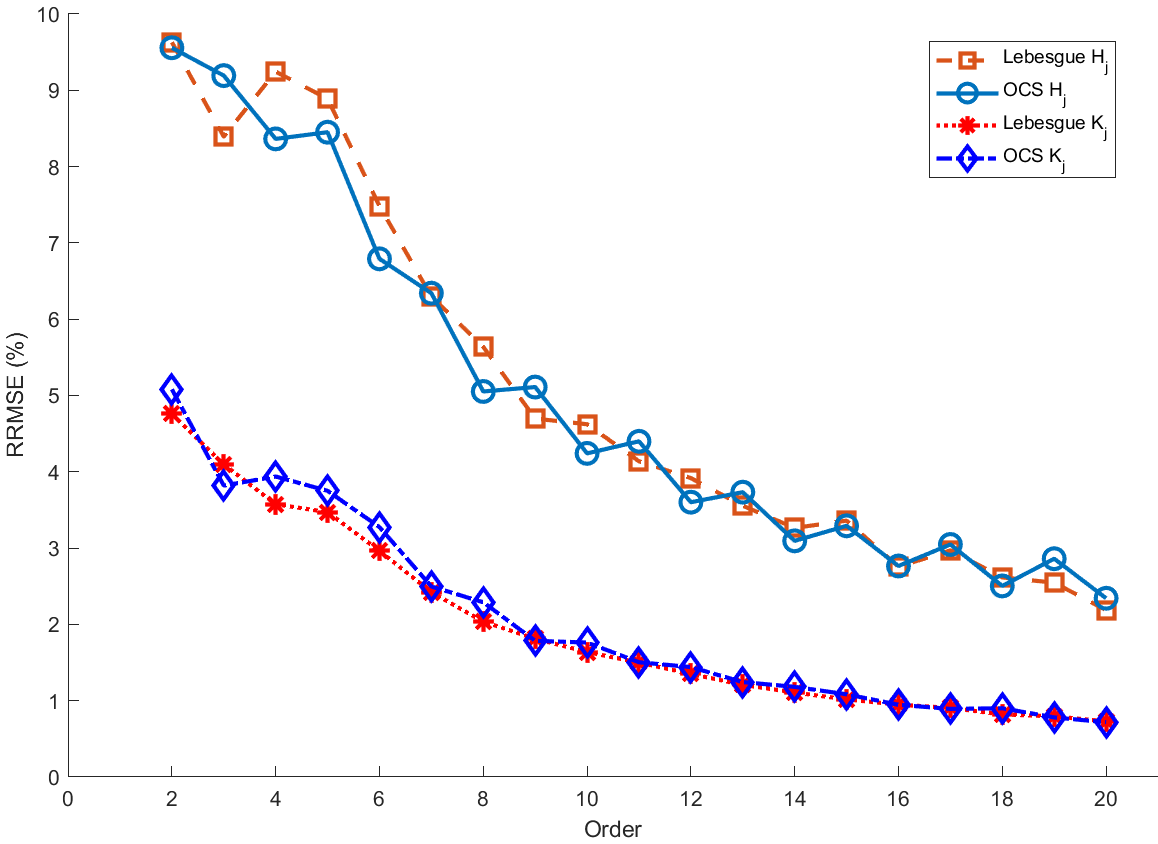}
    \includegraphics[width=6.7cm]{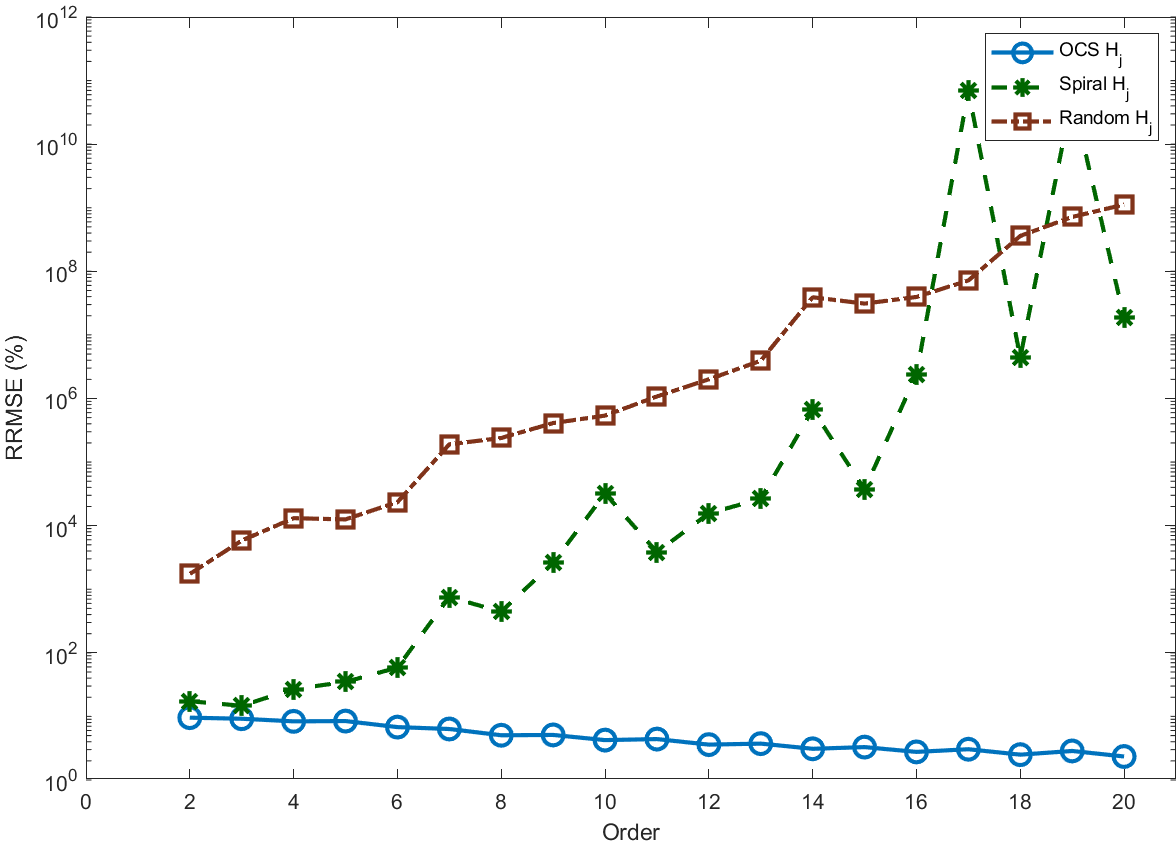}
    \caption{Relative root mean square error $(\%)$, for $H_j$ and $K_j$ and \textit{Lebesgue} and OCS samplings (left), and for $H_j$ and OCS, spiral and random samplings (right), with orders $n=2,\dotsc,20$.}
    \label{figure:rmsehj}
\end{figure}

Using these point configurations and the corresponding function values, we reconstruct the wavefront by solving the associated interpolation problem. The experiment is conducted for both orthogonal bases $H_j$ and $K_j$. For each interpolation order $n = 1, \dots, 30$ (which corresponds to a maximum of 496 nodes), we reconstruct 100 wavefront functions, each randomly generated according to equation~\eqref{eq:kolmo}. The relative root mean square error (RRMSE), 
\begin{equation}
    \text{RRMSE}:=\left(\frac{ \sum_{j=1}^M|\widetilde{f}^{approx}(\rho_j,\theta_j)-\widetilde{f}(\rho_j,\theta_j)|^2}{\sum_{j=1}^M |\widetilde{f}(\rho_j,\theta_j)|^2}\right)^{1/2},
\end{equation}
is evaluated on a uniform grid of $M\approx 2500$ points within each of the 36 hexagons which compose $\mathcal{H}$. Finally, we compute the average RRMSE for all 100 functions to obtain the mean interpolation error for each order $n$.

As the results indicate, we achieve a relative error of approximately $0.5\%$ when using the basis functions $K_j$ with both OCS and \textit{Lebesgue} sampling schemes, representing a highly accurate reconstruction. For the basis $H_j$, the error reaches approximately $1.6\%$, which still demonstrates a level of precision that is suitable for practical applications. The results for the transferred OCS and \textit{Lebesgue points}, using both orthogonal bases, are shown in Figure~\ref{figure:rmsehj} (left), for interpolation degrees up to n = 20.

In Figure~\ref{figure:rmsehj} (right), we also observe that both the translated spiral configuration and the random sampling method exhibit numerical instability, leading to significantly poorer performance compared to the translated OCS and \textit{Lebesgue points}. The resulting interpolation error increases rapidly with the degree $n$, indicating that these sampling strategies may lead to substantial approximation errors also in other practical applications.

\section{Conclusions}
\label{sec:conclusion}

In this work, we developed and analyzed sampling strategies for non-circular geometries\textendash specifically, the hexagon and the annulus\textendash tailored to Zernike-like orthogonal functions constructed via diffeomorphic mappings. The simpler case of the ellipse has also been studied. We provided theoretical results on the condition numbers of the associated collocation matrices and confirmed them through extensive numerical experiments. The proposed sampling patterns, when transferred from the unit disk, exhibit excellent numerical stability and enable high-order interpolation with low approximation error in the other domains.

A key application demonstrated the effectiveness of these methods in wavefront reconstruction for segmented mirror telescopes composed of hexagonal facets. Using realistic, Kolmogorov-generated wavefronts, we achieved accurate reconstructions with relative errors as low as 0.5\%, even at moderate polynomial degrees. These results highlight the practical viability of the approach in optical engineering and related fields.

Our results show that appropriately transferred sampling configurations, especially \textit{Optimal Concentric Sampling} and \textit{Lebesgue points}, allow stable high-order interpolation and effective wavefront modeling in complex optical systems. Moreover, the \textit{Optimal Concentric Samplings} can be computed with an explicit expression, which is a significant advantage in practice, whereas the \textit{Lebesgue points} need to be approximated numerically for every order.

\section*{Acknowledgements}
The first author (SDE) is supported by University of Almeria’s programme for research and knowledge transfer (PPIT-UAL).

The second author (AMF) was partially supported by Simons Foundation Collaboration Grants for Mathematicians (grant MPS-TSM-00710499).
He also acknowledges the support of the project PID2021-124472NB-I00, funded by MICIU/AEI/10.13039/501100011033 and by ``ERDF A way of making Europe'', as well as the support of Junta de Andaluc\'{\i}a (Instituto Interuniversitario Carlos I de F\'{\i}sica Te\'orica y Computacional).

The work of the third author (DRL) is supported by the projects PID2021-124472NB-I00 and PID2022-139293NB-C31 funded by MICIU/AEI/10.13039/501100011033 and by ``ERDF A way of making Europe''.

The authors thank the support from the research group FQM-229 (Junta de Andalucía), and from the Center for Development and Transfer of Mathematical Research to Industry CDTIME (University of Almería).


\begin{thebibliography}{10}
\expandafter\ifx\csname url\endcsname\relax
  \def\url#1{\texttt{#1}}\fi
\expandafter\ifx\csname urlprefix\endcsname\relax\def\urlprefix{URL }\fi
\expandafter\ifx\csname href\endcsname\relax
  \def\href#1#2{#2} \def\path#1{#1}\fi

\bibitem{Zernike1934}
F.~Zernike, Beugungstheorie des schneidenver-fahrens und seiner verbesserten
  form, der phasenkontrastmethode, Physica 1 (1934) 689--704.

\bibitem{Dai2008}
G.~Dai, Wavefront {O}ptics for {V}ision {C}orrection, SPIE Press, 2008.

\bibitem{Martinez-Finkelshtein:2009kx}
A.~Mart\'{\i}nez-Finkelshtein, A.~M. Delgado, G.~M. Castro, A.~Zarzo, J.~L.
  Ali\'o, Comparative analysis of some modal reconstruction methods of the
  shape of the cornea from the corneal elevation data, Invest. Ophthalmol. Vis.
  Sci. 50~(12) (2009) 5639--5645.

\bibitem{Ares2006}
M.~Ares, S.~Royo, Comparison of cubic {B}-spline and {Zernike}-fitting
  techniques in complex wavefront reconstruction, Applied Optics 45 (2006)
  6945--6964.

\bibitem{Schneider2009}
M.~Schneider, D.~R. Iskander, M.~J. Collins, Modeling corneal surfaces with
  rational functions for high-speed videokeratoscopy data compression, IEEE
  Transactions on Biomedical Engineering 56 (2) (2009) 493--499.

\bibitem{MR3683670}
A.~Mart\'{i}nez-Finkelshtein, D.~Ramos-L\'{o}pez, D.~R. Iskander, Computation
  of 2{D} {F}ourier transforms and diffraction integrals using {G}aussian
  radial basis functions, Appl. Comput. Harmon. Anal. 43~(3) (2017) 424--448.

\bibitem{Mahajan1994}
V.~Mahajan, Zernike {A}nnular {P}olynomials and {O}ptical {A}berrations of
  {S}ystems with {A}nnular {P}upils, Applied Optics 33 (1994) 8125--8127.

\bibitem{Mahajan2006}
V.~Mahajan, G.~Dai, Orthonormal polynomials for hexagonal pupils, Optics
  Letters 31 (2006) 2462--2464.

\bibitem{Mahajan2006a}
V.~Mahajan, G.~Dai, Orthonormal {P}olynomials in {W}avefront {A}nalysis:
  {A}nalytical {S}olution, Frontiers in Optics (2006).

\bibitem{Navarro2014}
R.~Navarro, J.~L. López, J.~A. Díaz, E.~P. Sinusía, Generalization of
  {Z}ernike polynomials for regular portions of circles and ellipses, Optics
  Express 22~(18) (2014) 21263--21279.

\bibitem{Ferreira2015}
C.~Ferreira, J.~L. López, R.~Navarro, E.~Pérez~Sinusía, Zernike-like systems
  in polygons and polygonal facets, Applied Optics 54~(21) (2015) 6575--6583.

\bibitem{Navarro2009}
R.~Navarro, J.~Arines, R.~Rivera, Direct and inverse discrete {Z}ernike
  transform, Optics Express {\bf 17}(26) (2009) 24269--24281.

\bibitem{Cuyt2012}
A.~Cuyt, I.~Yaman, B.~Ibrahimoglu, B.~Benouahmane, Radial orthogonality and
  {L}ebesgue constants on the disk, Numerical Algorithms 61 (2012) 291--313.

\bibitem{RamosLopez2016}
D.~Ramos-López, M.~A. S{\'a}nchez-Granero, M.~Fern{\'a}ndez-Mart{\'\i}nez,
  A.~Mart{\'\i}nez-Finkelshtein, Optimal sampling patterns for zernike
  polynomials, Applied Mathematics and Computation 274 (2016) 247--257.

\bibitem{Carnicer2014}
J.~Carnicer, C.~Godés, Interpolation on the disk, Numerical Algorithms 66
  (2014) 1--16.

\bibitem{Meurant2019}
G.~Meurant, A.~Sommariva, On the computation of sets of points with low
  {L}ebesgue constant on the unit disk, Journal of Computational and Applied
  Mathematics 345 (2019) 388--404.

\bibitem{vanKooten2022}
M.~A.~M. van Kooten, S.~Ragland, R.~Jensen-Clem, Y.~Xin, J.-R. Delorme,
  J.~Kent~Wallace, On-sky reconstruction of {K}eck {P}rimary {M}irror {P}iston
  {O}ffsets using a {Z}ernike {W}avefront {S}ensor, The Astrophysical Journal
  932~(2) (2022) 109.
\newblock \href {https://doi.org/10.3847/1538-4357/ac6ba2}
  {\path{doi:10.3847/1538-4357/ac6ba2}}.

\bibitem{Ragland2022}
S.~Ragland, P.~Wizinowich, M.~van Kooten, M.~Bottom, B.~Calvin, M.~Fitzgerald,
  P.~Hinz, R.~Jensen-Clem, D.~Mawet, E.~Peretz, {Residual wavefront control of
  segmented mirror telescopes}, in: L.~Schreiber, D.~Schmidt, E.~Vernet (Eds.),
  Adaptive Optics Systems VIII, Vol. 12185, International Society for Optics
  and Photonics, SPIE, 2022, p. 121850Y.
\newblock \href {https://doi.org/10.1117/12.2630269}
  {\path{doi:10.1117/12.2630269}}.

\bibitem{Niu2022}
K.~Niu, C.~Tian, Zernike polynomials and their applications, Journal of Optics
  24 (2022).

\bibitem{Bos2010}
L.~Bos, S.~D. Marchi, A.~Sommariva, M.~Vianello, Computing {M}ultivariate
  {F}ekete and {L}eja {P}oints by {N}umerical {L}inear {A}lgebra, SIAM Journal
  on Numerical Analysis 48 (2010).

\bibitem{SommarivaSets}
A.~Sommariva, Good interpolation sets on the unit disk, available at
  \url{https://www.math.unipd.it/~alvise/sets.html} (2023).

\bibitem{Ibrahimoglu2016}
B.~Ibrahimoglu, Lebesgue functions and {L}ebesgue constants in polynomial
  interpolation, Journal of Inequalities and Applications 93 (2016).

\bibitem{Bojanov2003}
B.~Bojanov, Y.~Xu, On polynomial interpolation of two variables, Journal of
  Approximation Theory 120 (2003) 267--282.

\bibitem{Cinzano2020}
P.~Cinzano, F.~Falchi, Toward an atlas of the number of visible stars, Journal
  of Quantitative Spectroscopy and Radiative Transfer 253 (2020).
\newblock \href {https://doi.org/10.1016/j.jqsrt.2020.107059}
  {\path{doi:10.1016/j.jqsrt.2020.107059}}.

\bibitem{Dudt2020}
D.~W. Dudt, E.~Kolemen, Desc: A stellarator equilibrium solver, Physics of
  Plasmas 27~(10) (2020) 102513.

\bibitem{Díaz2023}
J.~A. Díaz, R.~Navarro, Geometric optical transfer function and pupil sampling
  patterns, Optik 279 (2023).
\newblock \href {https://doi.org/10.1016/j.ijleo.2023.170746}
  {\path{doi:10.1016/j.ijleo.2023.170746}}.

\bibitem{Roddier1990}
N.~Roddier, Atmospheric wavefront simulation using {Z}ernike polynomials,
  Optical Engineering 29 (1990) 1174--1180.




\end{thebibliography}
\end{document}